\documentclass[a4paper,12pt,leqno]{article}
\usepackage{amsmath,amssymb,amsfonts}

\makeatletter
\newbox\bk@bxb
\newbox\bk@bxa
\newif\if@bkcont
\newcount\bk@lcnt

\def\breakboxskip{2pt}
\def\breakboxparindent{1.8em}

\def\breakbox{\vskip\breakboxskip\relax
\setbox\bk@bxb\vbox\bgroup
\advance\linewidth -2\fboxrule
\hsize\linewidth\@parboxrestore
\parindent\breakboxparindent\relax}

\def\bk@split{%
\@tempdimb\ht\bk@bxb 
\advance\@tempdimb\dp\bk@bxb
\setbox\bk@bxa\vsplit\bk@bxb to\z@ 
\setbox\bk@bxa\vbox{\unvbox\bk@bxa}
\setbox\@tempboxa\vbox{\copy\bk@bxa\copy\bk@bxb}
\advance\@tempdimb-\ht\@tempboxa
\advance\@tempdimb-\dp\@tempboxa}

\def\bk@addfsepht{%
\setbox\bk@bxa\vbox{\vskip\fboxsep\box\bk@bxa}}

\def\bk@addskipht{%
\setbox\bk@bxa\vbox{\vskip\@tempdimb\box\bk@bxa}}

\def\bk@addfsepdp{%
\@tempdima\dp\bk@bxa
\advance\@tempdima\fboxsep
\dp\bk@bxa\@tempdima}

\def\bk@addskipdp{%
\@tempdima\dp\bk@bxa
\advance\@tempdima\@tempdimb
\dp\bk@bxa\@tempdima}

\def\bk@line{%
\hbox to \linewidth{%
\hskip-2\fboxsep\vrule \@width\fboxrule\hskip.5\fboxsep\vrule \@width\fboxrule\hskip1.5\fboxsep
\box\bk@bxa\hfil
}}%

\def\endbreakbox{\egroup
\ifhmode\par\fi{\noindent\bk@lcnt\@ne
\@bkconttrue\baselineskip\z@\lineskiplimit\z@
\lineskip\z@\vfuzz\maxdimen
\bk@split\bk@addfsepht\bk@addskipdp
\ifvoid\bk@bxb 
\def\bk@fstln{\bk@addfsepdp
\hskip-\parindent\vbox{\llap{\raisebox{-2ex}{\rule{1.5\fboxsep}{\fboxrule}\hskip.5\fboxsep}}\bk@line\llap{\rule{1.5\fboxsep}{\fboxrule}\hskip.5\fboxsep}}}

\else 
\def\bk@fstln{\vbox{\llap{\raisebox{-2ex}{\rule{1.5\fboxsep}{\fboxrule}\hskip.5\fboxsep}}\bk@line}\hfil%
\advance\bk@lcnt\@ne
\loop
\bk@split\bk@addskipdp\leavevmode
\ifvoid\bk@bxb 
\@bkcontfalse\bk@addfsepdp
\vtop{\bk@line\llap{\rule{2\fboxsep}{\fboxrule}}}%

\else 
\bk@line
\fi
\hfil\advance\bk@lcnt\@ne
\if@bkcont\repeat}%
\fi
\leavevmode\bk@fstln\par}\vskip\breakboxskip\relax}

\def\smp{\smallskip\par}

\def\pf{\noindent{\bf Proof~:}\ }

\def\findemo{~\leaders\hbox to 1em{\hss\  \hss}\hfill~\raisebox{.5ex}{\framebox[1ex]{}}\smp}

\def\mpn{\medskip\par\noindent}
\def\smpn{\smallskip\par\noindent}

\def\smp{\smallskip\par}
\def\smpn{\smallskip\par\noindent}

\def\mpoint{\;\;.}
\def\mvirg{\;\;,}

\def\Inf{{\rm Inf}}

\def\Inf{{\rm Inf}}

\def\Aut{{\rm Aut}}

\def\Irr{{\rm Irr}}

\def\op{^{op}}

\def\N{\mathbb{N}}

\newcommand{\romain}[1]{\uppercase\expandafter{\romannumeral #1}}

\newcommand{\flh}[2]{\mathop{\hbox to 12mm{\rightarrowfill}}_{\displaystyle #2}^{\displaystyle #1}\limits}
\newcommand{\sflh}[2]{\mathop{\hbox to 12mm{\rightarrowfill}}_{\scriptstyle #2}^{\scriptstyle #1}\limits}

\newcommand{\sou}[1]{\underline{#1}}
\newcommand{\sur}[1]{\,\overline{\! #1}}

\def\op{^{op}}

\newcommand{\carre}[8]{\begin{array}{ccc}
#1&\mathop{\hbox to 12mm{\rightarrowfill}}^{\displaystyle{#2}}\limits&#3\\
\llap{$\displaystyle{#4}$}\left\downarrow\vbox to 6mm{}\right. & & \left\downarrow\vbox to 6mm{}\right.\rlap{$\displaystyle{#5}$}\\
#6&\mathop{\hbox to 12mm{\rightarrowfill}}_{\displaystyle #7}\limits&#8\\
\end{array}}

\newcommand{\carrem}[8]{\begin{array}{ccc}
#1&\mathop{\hbox to 12mm{\rightarrowfill}}^{\displaystyle #2}\limits&#3\\
\llap{$\displaystyle #4$}\left\uparrow\vbox to 6mm{}\right. & & \left\uparrow\vbox to 6mm{}\right.\rlap{$\displaystyle #5$}\\
#6&\mathop{\hbox to 12mm{\rightarrowfill}}_{\displaystyle #7}\limits&#8\\
\end{array}}

\newenvironment{enonce}[1]{\pagebreak[2]\refstepcounter{subsection}\refstepcounter{prop}\smpn{{\bf \thesection.\arabic{prop}.\ \ #1:}}\begin{it} }{\end{it}\smp}
\newenvironment{enonce*}[1]{\pagebreak[2]\smpn{#1:}\begin{it} }{\end{it}\smp}
\newcommand{\result}[1]{\begin{enonce}{#1}}
\def\fresult{\end{enonce}}
\newcommand{\npar}{\smallskip\par\noindent\pagebreak[2]\refstepcounter{subsection}\refstepcounter{prop}{\bf \thesection.\arabic{prop}.\ \ }}

\newenvironment{mth}[1]{\begin{breakbox}\begin{enonce}{#1}}{\end{enonce}\end{breakbox}}
\newenvironment{mth*}[1]{\begin{breakbox}\begin{enonce*}{#1}}{\end{enonce*}\end{breakbox}}
\newenvironment{rem}[1]{\refstepcounter{subsection}\refstepcounter{prop} \mpn{{\bf \thesection.\arabic{prop}.}\ \ \bf#1:}}{\smp}

\makeatletter
\renewenvironment{enumerate}{\ifnum \@enumdepth >3 \@toodeep\else
      \advance\@enumdepth \@ne
      \edef\@enumctr{enum\romannumeral\the\@enumdepth}\list
      {\csname label\@enumctr\endcsname}{\setlength{\topsep}{1ex}\setlength{\itemsep}{0pt}\usecounter
        {\@enumctr}\def\makelabel##1{\hss\llap{##1}}}\fi}{\endlist}
\renewenvironment{itemize}{\ifnum \@itemdepth >3 \@toodeep\else \advance\@itemdepth \@ne
\edef\@itemitem{labelitem\romannumeral\the\@itemdepth}%
\list{\csname\@itemitem\endcsname}{\setlength{\topsep}{1ex}\setlength{\itemsep}{0pt}\def\makelabel##1{\hss\llap{##1}}}\fi}
{\endlist}
\makeatother
\makeatletter
\def\@sect#1#2#3#4#5#6[#7]#8{\ifnum #2>\c@secnumdepth
    \let\@svsec\@empty\else
    \refstepcounter{#1}\edef\@svsec{\csname the#1\endcsname .\hskip .5em}\fi
    \@tempskipa #5\relax
     \ifdim \@tempskipa>\z@
       \begingroup #6\relax
         \@hangfrom{\hskip #3\relax\@svsec}{\interlinepenalty \@M #8\par}%
       \endgroup
      \csname #1mark\endcsname{#7}\addcontentsline
        {toc}{#1}{\ifnum #2>\c@secnumdepth \else
                     \protect\numberline{\csname the#1\endcsname}\fi
                   #7}\else
       \def\@svsechd{#6\hskip #3\relax  
                  \@svsec #8\csname #1mark\endcsname
                     {#7}\addcontentsline
                          {toc}{#1}{\ifnum #2>\c@secnumdepth \else
                            \protect\numberline{\csname the#1\endcsname}\fi
                      #7}}\fi
    \@xsect{#5}}
\def\section{\@startsection {section}{1}{\z@}{-3.5ex plus-1ex minus
    -.2ex}{2.3ex plus.2ex}{\reset@font\Large\bf}}  

\makeatother
\renewenvironment{equation}{\refstepcounter{subsection}\refstepcounter{prop}$$}{\leqno{\bf (\theprop)}$$}

\def\mar[#1]{\ar@{-}[#1]|-{\object@{<}}}
\def\marb[#1]{\ar@{-}[#1]|{\object+{  }}}

\usepackage{xcolor}
\usepackage{setspace}
\def\Sup{{\rm Sup}}
\def\endpf{\findemo}
\def\Grm{{\rm Grm}}
\def\Blb{{\mathrm{\widehat{G}}}}
\def\Idown{I_{\downarrow}}
\def\join{\mathop{\vee}\limits}
\def\meet{\mathop{\wedge}\limits}
\newcommand{\joinb}[2]{\mathop{\vee}\limits_{{\scriptstyle #1}\atop {\scriptstyle #2}}}
\newcommand{\meetb}[2]{\mathop{\wedge}\limits_{{\scriptstyle #1}\atop {\scriptstyle #2}}}
\usepackage[all]{xy}

\newcommand{\msqcup}{\,\sqcup\,}
\newcommand{\mcap}{\,\cap\,}
\newcommand{\mcup}{\,\cup\,}
\def\Idown{I_{\downarrow}}
\def\mp{\medskip\par}
\reversemarginpar
\renewenvironment{equation}{\refstepcounter{subsection}\refstepcounter{prop}$$}{\leqno{(\theprop)}$$}
\def\pointplein{\makebox[0ex]{$\bullet$}}
\def\pointcreux{\makebox[0ex]{$\circ$}}
\def\pointcarre{\makebox[0ex]{${\scriptscriptstyle\square}$}}

\begin{document}
\centerline{\Large\bf Germs in a poset}\vspace{2ex}
\centerline{\bf Serge Bouc}\vspace{2ex}
\begin{center}
\begin{minipage}{12cm}
\begin{footnotesize}
\begin{spacing}{1.05}
{\bf Abstract:} Motivated by the theory of correspondence functors, we introduce the notion of {\em germ} in a finite poset, and the notion of {\em germ extension} of a poset. We show that any finite poset admits a largest germ extension called its {\em germ closure}. We say that a subset $U$ of a finite lattice $T$ is {\em germ extensible} in $T$ if the germ closure of $U$ naturally embeds in $T$. We show that any for any subset $S$ of a finite lattice $T$, there is a unique germ extensible subset $U$ of $T$ such that $U\subseteq S\subseteq \sur{G}(U)$, where $\sur{G}(U)\subseteq T$ is the embedding of the germ closure of $U$.\vspace{1ex}
\end{spacing}
\end{footnotesize}
\end{minipage}
\end{center}
\section{Introduction}
In a series of papers with Jacques Th\'evenaz (\cite{corfun-finiteness}, \cite{corfun-lattices}, \cite{corfun-simple-projective}, \cite{corfun-tensor}, \cite{corfun-simplicity}), we develop the theory of {\em correspondence functors} over a commutative ring $k$, i.e. linear representations over $k$ of the category of finite sets, where morphisms are correspondences instead of maps. In this theory, finite lattices and finite posets play a crucial role, at various places. \par
In particular, we show (\cite{corfun-finiteness}, Theorem 4.7) that the simple correspondence functors are parametrized by triples $(E,R,V)$, where $E$ is a finite set, $R$ is a partial order relation on $E$ - that is, $(E,R)$ is a finite poset - and $V$ is a simple $k\Aut(E,R)$-module. Moreover, the evaluation at a finite set $X$ of the simple functor $S_{E,R,V}$ parametrized by the triple $(E,R,V)$ can be completely described (\cite{corfun-simplicity}, Theorem~6.6 and Theorem~7.9). It follows (\cite{corfun-simplicity}, Theorem~8.2) that when $k$ is a field, the dimension of $S_{E,R,V}(X)$ is given by 
\begin{equation}\label{dimension}\dim_k S_{E,R,V}(X)=\frac{\dim_kV}{|\Aut(E,R)|}\sum_{i=0}^{|E|}(-1)^i\binom{|E|}{i}\big(|G|-i\big)^{|X|}\mpoint
\end{equation}
The main consequence of these results is a complete description of the simple modules for the algebra over $k$ of the monoid of all relations on~$X$ (\cite{corfun-simplicity}, Section~8).\par
Formula~\ref{dimension} is obtained by first choosing a finite lattice $T$ such that the poset $\Irr(T)$ of join-irreducible elements of $T$ is isomorphic to the opposite poset $(E,R\op)$, and then constructing a specific subset $G=G_T$ of $T$ (see (\ref{def GT}) for a precise definition of $G_T$), which appears in the right hand side. Now for a given $X$, the left hand side of~(\ref{dimension}) only depends on the poset $(E,R)$ and the simple $k\Aut(E,R)$-module $V$, whereas the right hand side depends in addition on the cardinality of the set $G$, which {\em a priori} depends on $T$, and not just on $(E,R)$. We have checked (\cite{corfun-simplicity}~Corollary~6.7) that $|G|$ indeed only depends on $(E,R)$. A natural question is then to ask if the {\em subposet} $G$ of $T$ only depends on $(E,R)$, up to isomorphism, and not really on $T$ itself.\par
One of the aims of the present paper is to answer this question. In fact, the main aim is to introduce various structural results on posets and lattices, which appear to be new. The first notion we introduce is the notion of {\em germ} of a finite poset. A germ of a finite poset $S$ is an element of $S$ with specific properties (Definition~\ref{def germ}). \par
A related notion is the following: When $U$ is a subset of $S$, the poset $S$ is called a {\em germ extension} of $U$ if any element of $S-U$ is a germ of $S$. The main result of the paper (Theorem~\ref{germ closure}) is that conversely, being given a finite poset $U$, there is a (explicitly defined) finite poset $G(U)$, containing $U$ as a full subposet, which is the largest germ extension of $U$, in the following sense: First $G(U)$ is a germ extension of $U$, and moreover, if $S$ is a finite poset containing $U$ as a full subposet, and such that $S$ is a germ extension of $U$, then there is a unique full poset embedding $S\to G(U)$ which restricts to the identity map of $U$. For this reason, the poset $G(U)$ will be called the {\em germ closure} of $U$. \par
This poset $G(U)$ can be viewed as a structural answer to the above question: In the case of a finite lattice $T$ with poset $(E,R)$ of join-irreducible elements, we show (Theorem~\ref{claim}) that the poset set $G$ identifies canonically with $G(E,R)$, and in particular, it only depends on the poset $(E,R)$.\par
In Section 3, we consider {\em germ extensible} subsets of a finite lattice. For any full subposet $U$ of a finite lattice $T$, the inclusion map $U\hookrightarrow T$ extends to a canonical map of posets $\nu:G(U)\to T$. We say that $U$ is germ extensible in $T$ if this map $\nu$ is injective, and in this case, we denote by $\sur{G}(U)\subseteq T$ its image. We give a characterization of germ extensible subsets of a lattice (Theorem~\ref{germ extensible}), and then show (Theorem~\ref{partition}) that for any subset $S$ of $T$, there exists a unique germ extensible subset $U$ of $T$ such that $U\subseteq S\subseteq \sur{G}(U)$. \par
In other words, the poset of subsets of $T$ is partitioned by the intervals $[U,\sur{G}(U)]$, where $U$ is a germ extensible subset of $T$. In a forthcoming paper, we will show how this rather surprising result yields a natural filtration of the correspondence functor $F_T$ associated to $T$ (\cite{corfun-lattices}, Definition~4.1), by fundamental functors indexed by germ extensible subsets of~$T$.
\par
The last section of the paper (Section 5) lists some examples of germs, germ closures, and germ extensible subsets of lattices.\section{Germs in a poset}
Throughout the paper, we use the symbol $\subseteq$ for inclusion of sets, and the symbol $\subset$ for proper inclusion. We denote by $\sqcup$ the disjoint union of sets.\par
If $(U,\leq)$ is a poset, and $u,v$ are elements of $U$, we set
\begin{align*}
[u,v]_U&=\{w\in U\mid u\leq w\leq v\},\; &[u,v[_U&=\{w\in U\mid u\leq w< v\},\\
]u,v]_U&=\{w\in U\mid u< w\leq v\},\; &]u,v[_U&=\{w\in U\mid u< w< v\},\\
\,]\,.\,,v]_U&=\{w\in U\mid  w\leq v\},\; &\,]\,.\,,v[_U&=\{w\in U\mid  w< v\},\\
[v,\,.\,[_U&=\{w\in U\mid  w\geq v\},\; &]v,\,.\,[_U&=\{w\in U\mid  w> v\}.
\end{align*}
When $V$ is a subset of a poset $U$, we denote by $\Sup_UV$ the least upper bound of $V$ in $U$, when it exists. Similarly, we denote by $\Inf_UV$ the greatest lower bound of $V$ in $U$, when it exists. For $u\in U$, when we write $u=\Sup_UV$ (resp. $u=\Inf_UV$), we mean that $\Sup_UV$ exists (resp. that $\Inf_UV$ exists) and is equal to $u$. 
\begin{mth}{Definition}\label{def germ} Let $(U,\leq)$ be a finite poset. A {\em germ} of $U$ is an element $u\in U$ such that there exists an element $v\geq u$ in $U$ for which the following properties hold:
\begin{enumerate}
\item $u=\Sup_U\,]\,.\,,u[_U$ and $v=\Inf_U]v,\,.\,[_U$.
\item $[u,\,.\,[_U=[u,v]_U\msqcup ]v,\,.\,[_U$ and $\,]\,.\,,v]_U=\,]\,.\,,u[_U\msqcup [u,v]_U$.
\item $[u,v]_U$ is totally ordered.
\end{enumerate}
\end{mth}

\begin{mth}{Lemma} Let $(U,\leq)$ be a finite poset, and $u$ be a germ of $U$. Then there exists a unique element $v\geq u$ in $U$ with the properties of Definition~\ref{def germ}.\vspace{-1ex}
\end{mth}
\pf Let $v'\geq u$ be another element of $U$ with the properties of the element $v$ of Definition~\ref{def germ}. Then $v'\in [u,\,.\,[_U=[u,v]_U\msqcup ]v,\,.\,[_U$. In particular $u\leq v'\leq v$ or $v'\geq v$. If $v'\ngeq v$, then $u\leq v'<v$, and $]v',\,.\,[_U=]v',v]\msqcup ]v,\,.\,[_U$. Moreover $]v',v]_U$ has a smallest element $w$, since $[u,v]_U$ is totally ordered. Then $w$ is also the smallest element of $]v',\,.\,[_U$, thus $w=\Inf_U]v',\,.\,[_U\neq v'$, contradicting Property 1 of $v'$ in Definition~\ref{def germ}. It follows that $v'\geq v$. Exchanging the roles of $v$ and $v'$ gives $v\geq v'$, thus $v=v'$.\endpf
\begin{mth}{Definition and Notation}\label{def cogerm} Let $(U,\leq)$ be a finite poset.
\begin{enumerate} 
\item Let $u$ be a germ of $U$. The unique element $v\geq u$ of $U$ with the properties of Definition~\ref{def germ} is called the {\em cogerm} of $u$ in $U$.
\item The set of germs of $U$ is denoted by $\Grm(U)$.
\end{enumerate}
\end{mth}
\begin{rem}{Example} \label{germe cogerme}Let $u$ be an element of $U$ such that $u=\Sup_U\,]\,.\,,u[_U$ and $u=\Inf_U]u,\,.\,[_U$. Then $u$ is a germ of $U$, and $u$ is its own cogerm in $U$.
\end{rem}
\begin{rem}{Remark}\label{racines op} Clearly, if $u$ is a germ of the poset $U$, then its cogerm $v$ in~$U$ is a germ in the opposite poset $U\op$, and the cogerm of $v$ in $U\op$ is $u$. The correspondence $u\leftrightarrow v$ is a bijection between the germs of $U$ and the germs of~$U\op$.
\end{rem}
\begin{mth}{Lemma}\label{lemme racines} Let $(U,\leq)$ be a finite poset, let $u$ and $u'$ be distinct germs of~$U$, with respective cogerms $v$ and $v'$ in $U$. 
\begin{enumerate}
\item If $u'>u$, then $v'\geq u'>v\geq u$.
\item If $u'\leq v$, then $u'\leq v'< u\leq v$.
\end{enumerate}
\end{mth}
\pf 1) If $u'>u$, then $u'\in ]u,\,.\,[_U=]u,v]_U\msqcup]v,\,.\,[_U$. Moreover if $u'\leq v$, then $\,]\,.\,,u'[_U=\,]\,.\,,u[_U\msqcup[u,u'[_U$, and $[u,u'[_U$ has a greatest element $w$ since $[u,v]_U$ is totally ordered. Then $w$ is also the greatest element of $\,]\,.\,,u'[_U$, hence $w=\Sup_U\,]\,.\,,u'[_U\neq u'$, contradicting the first property of $u'$ in Definition~\ref{def germ}. It follows that $u'>v$. \mpn
2) If $u'\leq v$, then either $u'\leq u$ or $u'\in[u,v]_U$. If $u'\leq u$, then $u'<u$ since $u'\neq u$, and then $u'\leq v'<u\leq v$ by Assertion 1 (exchanging the roles of $(u,v)$ and $(u',v')$). Now if $u'\in[u,v]_U$, then $u'>u$ since $u\neq u'$, and then $u'>v$ by Assertion 1, contradicting the assumption $u'\leq v$.\findemo
\vspace{1ex}
\pagebreak[3]
\begin{mth}{Definition and Notation} Let $(S,\leq)$ be a finite poset, and let $U$ be a subset of $S$. For $B\subseteq S$, set $U_{\leq B}=\{u\in U\mid u\leq b\;\forall b\in B\}$, and for $s\in S$, set $U_{\leq s}=U_{\leq\{s\}}=\,]\,.\,,s]_S\mcap U$.
\begin{enumerate}
\item $U$ is said to {\em detect} $S$ if
$$\forall s,t\in S,\;\;\;s\leq t\;\;\; \iff \;\;\;U_{\leq s}\subseteq U_{\leq t}\mpoint$$
\item $S$ is called a {\em germ extension} of $U$ if $S=U\mcup \Grm(S)$.\vspace{-1ex}
\end{enumerate}
\end{mth}
\begin{mth}{Lemma}\label{detecte} Let $(S,\leq)$ be a finite poset, and let $U$ be a subset of $S$. If $S$ is a germ extension of $U$, then $U$ detects $S$.
\end{mth}
\pf Clearly if $s,t\in S$ and $s\leq t$, then $U_{\leq s}\subseteq U_{\leq t}$. Conversely, for $t\in S$, set
$$\Omega_t=\{s\in S\mid \;U_{\leq s}\subseteq U_{\leq t}\;\hbox{but}\;s\nleq t\}\mpoint$$
If $\Omega_t\neq\emptyset$, let $s$ be a minimal element of $\Omega_t$. If $s\in U$, then $s\in U_{\leq s}$, thus $s\in U_{\leq t}$, hence $s\leq t$, contradicting the assumption $s\in\Omega_t$. Hence $s\notin U$, thus $s\in \Grm(S)$, and in particular $s=\Sup_S\,]\,.\,,s[_S$. Now for $x\in \,]\,.\,,s[_S$, we have $U_{\leq x}=\,]\,.\,,x]_S\mcap U\subseteq \,]\,.\,,s]_S\mcap U=U_{\leq s}\subseteq U_{\leq t}$, hence $x\leq t$ by minimality of $s$. It follows that $t\geq \Sup_S\,]\,.\,,s[_S=s$, a contradiction. Hence $\Omega_t=\emptyset$ for any $t\in S$, as was to be shown.\endpf
\begin{mth}{Theorem}\label{S dans G} Let $(S,\leq)$ be a finite poset, and let $U$ be a full subposet of~$S$ such that $S$ is a germ extension of $U$. Let $s\in S$. Then one and only one of the following holds:
\begin{enumerate}
\item there exists a subset $B$ of $U$ such that $U_{\leq s}=U_{\leq B}$.
\item there exists a germ $r$ of $U$ such that $U_{\leq s}=\,]\,.\,,r[_U$.
\end{enumerate}
\end{mth}
\pf Let us first prove that one of the assertions 1 or 2 holds. If $s\in U$, then $U_{\leq s}=\,]\,.\,,s]_U$, so 1 holds for $B=\{s\}$. Now if $s\notin U$, then $s\in\Grm(S)$. Let $\hat{s}$ be the cogerm of $s$ in $S$. The totally ordered poset $[s,\hat{s}]$ is of the form
$$[s,\hat{s}]_S=\{s=s_0<s_1<\ldots<s_n=\hat{s}\}\mpoint$$
Consider first an element $u\in U$ such that $u\leq u'$ for any $u'\in ]\hat{s},\,.\,[_S\mcap U$. Suppose that there exists $t\in S$ such that $\hat{s}<t$ but $u\nleq t$, and let $t$ be maximal with this property. Then in particular $t\notin U$. Hence $t\in \Grm(S)$. Moreover $u\leq x$ for any $x\in ]t,\,.\,[_S$, by maximality of $t$. Let $\hat{t}$ be the cogerm of $t$ in $S$.\par
If $t=\hat{t}$, then $u\leq \Inf_S]t,\,.\,[_S=t$, a contradiction. Thus $t<\hat{t}$, and $[t,\hat{t}]=\{t=t_0<t_1<\ldots<t_m=\hat{t}\}$, for some $m\geq 1$. Then $t_i\in U$ for $i\geq 1$, for otherwise $t_i$ is a germ of $S$, and $t_i>t=t_0$, thus $t_i>\hat{t}=t_m$ by Lemma~\ref{lemme racines}, a contradiction. In particular $u\leq t_1$. If $u<t_1$, then $u\leq t_0=t$, a contradiction. Thus $u=t_1$. It follows that $t_1$ is the smallest element of $]\hat{s},\,.\,[_S\mcap U$. Since $t$ is the greatest element of $\,]\,.\,,t_1[_S$, it follows that $t$ is unique. In other words the set $\{x\in S\mid x>\hat{s},\; u\nleq x\}$ has a unique maximal element, i.e. it has a greatest element $t<t_1=u$. Thus
$$]\hat{s},\,.\,[_S=]\hat{s},t]_S\msqcup[u,\,.\,[_S\mpoint$$
Suppose that $]\hat{s},t[_S\neq\emptyset$, and let $t'$ be a maximal element of that poset. Then $t'\notin U$, for if $t'\in U$, then $t'\geq u>t$, since $t'>\hat{s}$. Thus $t'\in \Grm(S)$, and $t'<t$. By Lemma~\ref{lemme racines}, it follows that $t'\leq\hat{t}'<t$, where $\hat{t}'$ is the cogerm of $t'$ in $S$. Hence $t'=\hat{t'}$. But $]t',\,.\,[_S=]t',t]_S\msqcup[u,\,.\,[_S=\{t\}\msqcup [u,\,.\,[_S$, and $t=\Inf_S]t',\,.\,[_S$, a contradiction since $t'=\hat{t'}=\Inf_S]t',\,.\,[_S$.\par
Thus $]\hat{s},t[_S=\emptyset$. It follows that $]\hat{s},\,.\,[_S=\{t\}\msqcup[u,\,.\,[_S$, hence $\Inf]\hat{s},\,.\,[_S=t\neq\hat{s}$, a contradiction.\par
Therefore $u\leq t$ for any $t\in]\hat{s},\,.\,[_S$, hence $u\leq \hat{s}=\Inf_S]\hat{s},\,.\,[_S$. This shows that for $u\in U$
\begin{equation}\label{eq4}
u\leq u'\;\;\forall u'\in]\hat{s},\,.\,[_S\mcap U \iff u\leq\hat{s}\mpoint
\end{equation}

Now there are two cases:
\begin{itemize}
\item if $s=\hat{s}$, then by~(\ref{eq4})
$$u\leq s=\hat{s} \iff u\leq u'\;\;\forall u'\in]s,\,.\,[_S\mcap U\mvirg$$
so we are in Case 1 of the theorem, for $B=]s,\,.\,[_S\mcap U$.
\item if $s<\hat{s}$, that is $n\geq 1$. Then for $i\in\{1,\ldots,n\}$, the element $s_i$ is in $U$: indeed if $s_i\notin U$ then $s_i\in\Grm(S)$, and $s_i>s=s_0$. Hence $s_i>\hat{s}=s_n$ by Lemma~\ref{lemme racines}, a contradiction.\par
Now if $u\in U$, and $u\geq u'$ for any $u'\in \,]\,.\,,s_1[_U$, then $u\geq u'$ for any $u'\in \,]\,.\,,s_0]_S\mcap U$, that is $U_{\leq s_0}\subseteq U_{\leq u}$, hence $s_0\leq u$ by Lemma~\ref{detecte}. Thus $s_0<u$, i.e. $s_1\leq u$, showing that 
\begin{equation}\label{eq1}
s_1=\Sup_U\,]\,.\,,s_1[_U\mpoint
\end{equation}
Moreover 
\begin{equation}\label{eq2}
[s_1,\,.\,[_U=[s_0,\,.\,[_S\mcap U=\big([s_0,\hat{s}]_S\mcap U\big)\msqcup]\hat{s},\,.\,[_U=[s_1,\hat{s}]_U\msqcup]\hat{s},\,.\,[_U\mpoint
\end{equation}
Similarly
\begin{equation}\label{eq3}
\,]\,.\,,\hat{s}]_U=\big(\,]\,.\,,s_0]_S\mcap U\big)\msqcup \big([s_0,\hat{s}]_S\mcap U\big)=\,]\,.\,,s_1[_U\msqcup[s_1,\hat{s}]_U\mpoint
\end{equation}
Finally, from~(\ref{eq4}), we have that 
\begin{equation}\label{eq5}
\hat{s}=\Inf_U]\hat{s},\,.\,[_U\mpoint
\end{equation}
Now it follows from~(\ref{eq1}),~(\ref{eq2}),~(\ref{eq3}) and~(\ref{eq5}) that $s_1$ is a germ of~$U$, with cogerm $\hat{s}$. Moreover $U_{\leq s}=\{u\in U\mid u<s_1\}$, so this is Case~2 of the theorem.
\end{itemize}
It remains to see that the two cases of the theorem cannot occur simultaneously. So suppose that there exists $s\in S$, $B\subseteq U$, and $r\in\Grm(U)$ such that for any $u\in U$
$$u\leq s \iff u\leq b,\;\forall b\in B \iff u<r\mpoint$$
Then for $b\in B$, we have $b\geq u$ for any $u\in \,]\,.\,,r[_U$, hence $b\geq \Sup_U\,]\,.\,,r[_U=r$. Thus $r\leq b$ for any $b\in B$, hence $r<r$, a contradiction. This completes the proof of Theorem~\ref{S dans G}.\endpf
\begin{mth}{Theorem}\label{reciproque} Let $(S,\leq)$ be a finite poset, and $U$ be a full subposet of~$S$. Suppose that $U$ detects $S$, and that for any $s\in S$, one of the following holds:
\begin{enumerate}
\item there exists a subset $B$ of $U$ such that $U_{\leq s}=U_{\leq B}$.
\item there exists a germ $r$ of $U$ such that $U_{\leq s}=\,]\,.\,,r[_U$.
\end{enumerate}
Then $S$ is a germ extension of $U$.
\end{mth}
\pf Let $s\in S-U$. \mpn
$\bullet$ \sou{Step 1}: Let $t\in S$ such that $s'\leq t$ for all $s'\in\,]\,.\,,s[_S$. If $u\in U_{\leq s}$, that is if $u\in U$ and $u\leq s$, then $u<s$, and then $u\leq t$. In other words $U_{\leq s}\subseteq U_{\leq t}$, thus $s\leq t$, since $U$ detects $S$. This shows that $s=\Sup_S\,]\,.\,,s[_S$.\mpn
$\bullet$ \sou{Step 2}: Suppose first that there exists $B\subseteq U$ such that $U_{\leq s}=U_{\leq B}$. Let $t\in S$ such that $t\leq s'$, for all $s'\in]s,\,.\,[_S$. If $b\in B$, then $U_{\leq s}\subseteq U_{\leq b}$, hence $s\leq b$, since $U$ detects $S$, and $s<b$ since $s\notin U$. It follows that $t\leq b$ for all $b\in B$. Hence $U_{\leq t}\subseteq U_{\leq s}$, so $t\leq s$. This shows that $s=\Inf_S]s,\,.\,[_S$. Since $s=\Sup_S\,]\,.\,,s[_S$, it follows (see Example~\ref{germe cogerme}) that $s$ is a germ of $S$, equal to its cogerm in $S$.\mpn
$\bullet$ \sou{Step 3}: Suppose now that there exists a germ $r$ of $U$ such that $U_{\leq s}=\,]\,.\,,r[_U$. In particular $U_{\leq s}\subset \,]\,.\,,r]_U=U_{\leq r}$, so $s\leq r$, since $U$ detects $S$, and $s<r$ since $s\notin U$ and $r\in U$.\par
Let $t\in ]s,\,.\,[_S$. Let us show that $r\leq t$. Suppose first that there exists $B\subseteq U$ such that $U_{\leq t}=U_{\leq B}$. Let $u\in U$ with $u<r$, and let $b\in B$. Then $u\leq s$, hence $u<t$, thus $u\leq b$. It follows that $b\geq \Sup_U\,]\,.\,,r[_U=r$. This holds for any $b\in B$, so $r\leq t$.\par
Assume now that there is a germ $r'$ of $U$ such that $U_{\leq t}=\,]\,.\,,r'[_U$. Then $U_{\leq s}$ is strictly contained in $U_{\leq t}$ since $s< t$ and $U$ detects $S$. Thus $\,]\,.\,,r[_U\subset \,]\,.\,,r'[_U$. In particular $r'\geq u$ for any $u\in \,]\,.\,,r[_U$, hence $r'\geq \Sup_U\,]\,.\,,r[_U=r$, and $r'>r$ since $\,]\,.\,,r[_U\subset \,]\,.\,,r'[_U$. Thus $r\in\,]\,.\,,r'[_U=U_{\leq t}$, i.e. $r\leq t$ again.\par
This shows that $r$ is the smallest element of $]s,\,.\,[_S$. Let $\hat{r}$ be the cogerm of $r$ in $U$. Then $[r,\hat{r}]_U=\{r=r_0<r_1<\ldots<r_n=\hat{r}\}$ (possibly $n=0$ if $r=\hat{r}$). \mpn
$\bullet$ \sou{Step 4}: Let $t\in S$ with $t>s$, i.e. $t\geq r$, and suppose that $\hat{r}\nleq t$. Then there exists an integer $m\in \{0,\ldots n-1\}$ such that $r_m\leq t$ but $r_{m+1}\nleq t$.\par
Suppose first that there exists $B\subseteq U$ such that $U_{\leq t}=U_{\leq B}$. Thus $r_m\leq b$ for all $b\in B$, but there exists $b_0\in B$ such that $r_{m+1}\nleq b_0$. Hence $B\subseteq [r_m,\,.\,[_U=[r_m,\hat{r}]_U\msqcup]\hat{r},\,.\,[_U$. If $r_m\notin B$, then $B\subseteq [r_{m+1},\hat{r}]_U\msqcup]\hat{r},\,.\,[_U$, as $r_{m+1}$ is the smallest element of $]r_m,\,.\,[_U$. This contradicts $r_{m+1}\nleq b_0$. Hence $r_m\in B$, and
$$U_{\leq t}\subseteq\{u\in U\mid u\leq r_m\}=U_{\leq r_m}\mpoint$$
It follows that $t\leq r_m$, since $U$ detects $S$. Therefore $t=r_m$ in this case, since we had $r_m\leq t$.\par
Suppose now that there is a germ $r'$ of $U$ such that $U_{\leq t}=\,]\,.\,,r'[_U$. Then $r<r'$, since $r\leq t$. Hence $\hat{r}<r'$ by Lemma~\ref{lemme racines}. Then $U_{\leq \hat{r}}\subseteq \,]\,.\,,r'[_U=U_{\leq t}$, and $\hat{r}\leq t$ since $U$ detects $S$. This contradicts the assumption on $t$.\par
This shows that if $t\in S$ and $t\geq s$, then $t\in\{s<r=r_0,\ldots,r_n=\hat{r}\}$ or $t> \hat{r}$.\mpn
$\bullet$ \sou{Step 5}: Let $t\in S$ with $t\leq\hat{r}$. Then $U_{\leq t}\subseteq \,]\,.\,,\hat{r}]_U=\,]\,.\,,r[_U\msqcup [r,\hat{r}]_U$. If $U_{\leq t}\subseteq \,]\,.\,,r[_U=U_{\leq s}$, then $t\leq s$ since $U$ detects $S$. Otherwise $U_{\leq t}$ has a greatest element $r_m\in\{r=r_0,r_1,\ldots,r_n=\hat{r}\}$, and $U_{\leq t}=U_{\leq r_m}$, thus $t=r_m$.\par
Hence if $t\in S$ and $t\leq \hat {r}$, then $t<s$ or $t\in\{s<r=r_0,\ldots,r_n=\hat{r}\}$.\mpn
$\bullet$ \sou{Step 6}: Let $t\in S$ such that $t\leq u$ for all $u\in U$ with $u>\hat{r}$. If $u'\in U_{\leq t}$, then $u'\leq u$ for all $u\in U$ with $u>\hat{r}$, hence $u'\leq\Inf_U]\hat{r},\,.\,[_U=\hat{r}$. It follows that $U_{\leq t}\subseteq \,]\,.\,,\hat{r}]_U=U_{\leq\hat{r}}$, hence $t\leq \hat{r}$, as $U$ detects $S$.\par
This shows {\em a fortiori} that if $t\in S$ and $t\leq s'$ for all $s'\in]\hat{r},\,.\,[_S$, then $t\leq\hat{r}$, that is $\hat{r}=\Inf_S]\hat{r},\,.\,[_S$.\mpn
$\bullet$ \sou{Step 7}: The conclusion of Steps 1, 4, 5 and 6 above show that $s$ is a germ of~$S$ under the assumption of Step 3, with cogerm $\hat{r}$ in $S$. Together with Step~2, this shows that $S-U\subseteq \Grm(S)$, in other words, that $S$ is a germ extension of $U$.\findemo
\begin{mth}{Corollary}\label{sous extension} Let $(S,\leq)$ be a finite poset, and $U\subseteq S$. If $S$ is a germ extension of $U$, so is any full subposet $R$ of $S$ containing $U$.
\end{mth}
\pf Let $R$ be a full subposet of $S$ containing $U$. Then $U$ detects $S$ by Lemma~\ref{detecte}, so $U$ detects $R\subseteq S$. Moreover, by Theorem~\ref{S dans G}, for any $x\in R$, the set $U_{\leq x}$ is equal to $\,]\,.\,,r[_U$ for some germ $r$ of the (full) subposet $U$ of $S$ (which is also a full subposet of $R$), or there exists a subset $B$ of $U$ such that $U_{\leq x}=U_{\leq B}$. By Theorem~\ref{reciproque}, it follows that $R$ is a germ extension of $U$.\findemo\pagebreak[3]
\begin{mth}{Definition and Notation}\label{def G} Let $(U,\leq)$ be a finite poset. Set
\begin{eqnarray*}
\Lambda(U)&=&\{s\subseteq U\mid\exists B\subseteq U,\;s=U_{\leq B}\},\\\Blb(U)&=&\{s\subseteq U\mid\exists r\in\Grm(U),\;s=\,]\,.\,,r[_U\}\mpoint
\end{eqnarray*}
Let $G(U)=\Lambda(U)\mcup \Blb(U)$, considered as a full subposet of the poset $\Idown(U)$ of lower-subsets of $U$ (ordered by inclusion of subsets of $U$). The poset $G(U)$  is called the {\em germ closure} of $U$.
\end{mth}
\begin{rem}{Remark} \label{defbulb}The sets $\Lambda(U)$ and $\Blb(U)$ are special cases, in the case of the lattice $\Idown(U)$ of lower-subsets of $U$, of constructions one can define in an arbitrary finite lattice. The set $\Lambda(U)$ is the set of intersections (i.e. meets) of lower intervals of $U$, i.e. join-irreducible elements of $\Idown(U)$, since $U_{\leq B}=\mathop{\bigwedge}_{b\in B}\limits\,]\,.\,,b]_U$ (see Section~\ref{germs and lattices} for details).
\end{rem} 
The terminology {\em germ closure} is motivated by the following:
\begin{mth}{Theorem}\label{germ closure} Let $(U,\leq)$ be a finite poset.
\begin{enumerate}
\item The map $u\in U \mapsto \sou{u}=\,]\,.\,,u]_U$ is an isomorphism from $U$ onto a full subposet $\sou{U}$ of $G(U)$.
\item If $S$ is a full subposet of $G(U)$ containing $\sou{U}$, then $S$ is a germ extension of $\sou{U}$. In particular $G(U)$ is a germ extension of $\sou{U}$.
\item Let $(S,\leq)$ be a poset containing $U$ as a full subposet. If $S$ is a germ extension of $U$, then there exists a unique isomorphism of posets $j:S\to S'$ onto a full subposet $S'$ of $G(U)$ such that $j(u)=\sou{u}$ for all $u\in U$.
\end{enumerate}
\end{mth}
\pf For Assertion 1, observe that $\,]\,.\,,u]_U=\{u'\in U\mid u'\leq u\}$, so $\,]\,.\,,u]_U$ indeed belongs to $G(U)$. Moreover for $u,v\in U$, the inclusion $\,]\,.\,,u]_U\subseteq \,]\,.\,,v]_U$ is equivalent to $u\leq v$.\mpn
For Assertion 2, by Corollary~\ref{sous extension}, it suffices to consider the case $S=G(U)$. First it is clear that $\sou{U}$ detects $G(U)$: indeed if $s\in G(U)$ and $u\in U$, then $\sou{u}\subseteq s$ if and only if $u\in s$, because $s$ is a lower subset. In other words $\sou{U}_{\leq s}=\{\sou{u}\mid u\in s\}$.\par
If there exists a subset $B$ of $U$ such that $s=U_{\leq B}$, then 
$$\sou{U}_{\leq s}=\{\sou{u}\mid u\leq b\;\;\forall b\in B\}=\{\sou{u}\in\sou{U}\mid \sou{u}\subseteq \sou{b}\;\;\forall \sou{b}\in \sou{B}\}\mvirg$$
where $\sou{B}=\{\sou{b}\mid b\in B\}\subseteq \sou{U}$.\par
Otherwise there exists a germ $r$ of $U$ such that $s=\,]\,.\,,r[_U$. Clearly in this case 
$$\sou{U}_{\leq s}=\{\sou{u}\mid u<r\}=\{\sou{u}\in\sou{U}\mid \sou{u}\subset\sou{r}\}\mvirg$$
and $\sou{r}$ is a germ of $\sou{U}$ since $u\mapsto \sou{u}$ is an isomorphism  of posets from $U$ to $\sou{U}$.\par
Now the assumptions of Theorem~\ref{reciproque} are fulfilled. It follows that $G(U)$ is a germ extension of $\sou{U}$, as was to be shown.\mpn
For Assertion 3, let $(S,\leq)$ be a finite poset containing $U$ as a full subposet, and assume that $S$ is a germ extension of $U$. Then $U$ detects $S$ by Lemma~\ref{detecte}, and moreover, for any $s\in S$, the set $U_{\leq s}$ belongs to $G(U)$, by Theorem~\ref{S dans G}.\par
For $s\in S$, set $j(s)=U_{\leq s}\in G(U)$, and $S'=j(S)$. Then $S'$ is a full subposet of $G(U)$ (because $U$ detects $S$), and $j$ is an isomorphism of posets $S\to S'$. Moreover $j(u)=U_{\leq u}=\sou{u}$ for $u\in U$, so $S'$ contains $\sou{U}$. This shows the existence of $S'$ and $j$ in Assertion~3. \par
For the uniqueness, let $j':S\to S''$ be an isomorphism of posets from $S$ to a full subposet $S''$ of $G(U)$, such that $j'(u)=\sou{u}$ for all $u\in U$. Then for $s\in S$ and $u\in U$
$$\sou{u}=j'(u)\subseteq j'(s)\iff u\leq s\iff j(u)=\sou{u}\subseteq j(s)\mpoint$$
In other words $\sou{U}_{\leq j'(s)}=\sou{U}_{\leq j(s)}$. Since $\sou{U}$ detects $G(U)$ by Assertion~2 and Lemma~\ref{detecte}, it follows that $j(s)=j'(s)$, hence $j=j'$ and $S'=S''$.\endpf
\begin{mth}{Proposition} Let $(U,\leq)$ be a finite poset. \begin{enumerate}
\item The poset $G(U)$ is the disjoint union of $\Lambda(U)$ and $\Blb(U)$.
\item If $s,t\in G(U)$, then $s\mcap t\in G(U)$. More precisely, if $s\nsubseteq t$ and $t\nsubseteq s$, then $s\mcap t\in \Lambda(U)$.
\item The subsets $\emptyset$ and $U$ belong to $G(U)$. Thus $G(U)$ is a lattice: for $s,t\in G(U)$, the infimum $s\meet t$ of $\{s,t\}$ in $G(U)$ is $s\mcap t$, and the supremum $s\join t$ is the intersection of all $x\in G(U)$ such that $x\supseteq s\mcup t$.
\end{enumerate}
\end{mth}
\pf By Definition~\ref{def G}, the poset $G(U)$ is the union of $\Lambda(U)$ and $\Blb(U)$, and by Theorem~\ref{S dans G}, this union is disjoint. Assertion~1 follows.\mpn
If $A$ and $B$ are subsets of $U$, then $U_{\leq A}\mcap U_{\leq B}=U_{\leq(A\mcup B)}$. Moreover, if $r$ is a germ of $U$, and if $r\in U_{\leq A}$, then $\,]\,.\,,r[_U\subset U_{\leq A}$. And if $r\notin U_{\leq A}$, then $U_{\leq A}\mcap \,]\,.\,,r[_U=U_{\leq A}\mcap \,]\,.\,,r]_U=U_{\leq(A\mcup\{r\})}$. Finally, if $r'$ is another germ of~$U$, and if $r\leq r'$, then $\,]\,.\,,r[_U\subseteq \,]\,.\,,r'[_U$. And if $r\nleq r'$ and $r'\nleq r$, then $\,]\,.\,,r[_U\mcap \,]\,.\,,r'[_U=\,]\,.\,,r]_U\mcap\,]\,.\,,r']_U=U_{\leq\{r,r'\}}$. Assertion~2 follows.\mpn
The set $U$ is equal to $U_{\leq\emptyset}$, so $U\in G(U)$. If $U$ has no smallest element, then $\emptyset=U_{\leq U}$, so $\emptyset\in G(U)$. And if $U$ has a smallest element $u$, then set $u_0=u$, and define inductively the finite sequence $u_0<u_1<\ldots<u_n$, where each $u_i$ is the smallest element of $]u_{i-1},\,.\,[_U$, for $i\geq 1$, and $]u_n,\,.\,[_U$ has no smallest element. Then $u=\Sup_U\,]\,.\,,u[_U$, since $\Sup_U(\emptyset)$ is the smallest element of~$U$. Moreover if $x\in U$, then either $x=u_i$ for some $i$, or $x> u_n$. And if $x\leq y$ for all $y\in ]u_n,\,.\,[_U$, then $x\notin ]u_n,\,.\,[_U$, since $]u_n,\,.\,[_U$ has no smallest element. Thus $x\leq u_n$, hence $u_n=\Inf_U]u_n,\,.\,[_U$. This shows that $u$ is a germ of~$U$, with cogerm $u_n$. Then $\,]\,.\,,u[_U=\emptyset\in G(U)$ in this case also. The last part of Assertion~3 now follows from Assertion~2, since $G(U)$ has a greatest element~$U$.\endpf
\begin{mth}{Theorem}\label{extension racines} Let $(S,\leq)$ be a finite poset, let $U$ be a full subposet of~$S$, and assume that $S$ is a germ extension of $U$. Then:
\begin{enumerate}
\item If $s$ is a germ of $S$ and if $s\in U$, then $s$ is a germ of $U$. In other words $U\mcap\Grm(S)\subseteq \Grm(U)$.
\item Let $r$ be a germ of $U$, with cogerm $\hat{r}\in U$. Then $[r,\hat{r}]_S=[r,\hat{r}]_U$, and one of the following holds:
\begin{enumerate}
\item $r=\Sup_S\,]\,.\,,r[_S$, and then $r$ is a germ of $S$, with cogerm $\hat{r}$ in $S$.
\item $\,]\,.\,,r[_S$ has a greatest element $s$. Then $r$ is the smallest element of $]s,\,.\,[_S$, and $s\in S-U$ is a germ of $S$, with cogerm $\hat{r}$ in $S$. In particular $r$ is not a germ of $S$.
\end{enumerate}
\end{enumerate}
\end{mth}
\pf  For Assertion~1, let $s\in U\mcap\Grm(S)$. If $u\in U$ is such that $u\geq v$ for any $v\in \,]\,.\,,s[_U$, let $t\in \,]\,.\,,s[_S$ be minimal such that $t\nleq u$. Then $t\notin U$, hence $t\in\Grm(S)$. In particular $t=\Sup_S\,]\,.\,,t[_S$. Moreover $t'\leq u$ for any $t'\in\,]\,.\,,t[_S$, by minimality of $t$. Hence $\Sup_S\,]\,.\,,t[_S=t\leq u$, a contradiction, which proves that $t\leq u$ for any $t\in \,]\,.\,,s[_S$. But then $\Sup_S\,]\,.\,,s[_S=s\leq u$. This shows that 
\begin{equation}\label{sup}s=\Sup_U\,]\,.\,,s[_U\mpoint
\end{equation}
Now let $\hat{s}$ be the cogerm of $s$ in $S$. There is an integer $n\in\N$ and elements $s_i\in S$, for $0\leq i\leq n$, such that
$${}[s,\hat{s}]_S=\{s=s_0<s_1<\ldots<s_n=\hat{s}\}\mpoint$$
If $s_i\notin U$ for some $i\in\{1,\ldots,n\}$, then $s_i\in\Grm(S)$, and $s_i>s$. By Lemma~\ref{lemme racines} $s_i>\hat{s}=s_n$, a contradiction. Thus $s_i\in U$ for any $i\in\{0,\ldots,,n\}$, i.e. $[s,\hat{s}]_S=[s,\hat{s}]_U$.\par
Now if $u\in U$ and $u\geq s$ then $u\in [s,\hat{s}]_S=[s,\hat{s}]_U$, or $u\geq\hat{s}$. Similarly if $u\leq\hat{s}$, then $u\in [s,\hat{s}]_S=[s,\hat{s}]_U$ or $u\leq s$. Thus
\begin{equation}\label{intervalle}
[s,\,.\,[_U=[s,\hat{s}]_U\msqcup]\hat{s},\,.\,[_U\;\;\hbox{and}\;\;\,]\,.\,,\hat{s}]_U=\,]\,.\,,s[_U\msqcup[s,\hat{s}]_U\mpoint
\end{equation}
Now let $u\in U$ such that $u\leq v$ for any $v\in]\hat{s},\,.\,[_U$. Suppose that there exists $t\in ]\hat{s},\,.\,[_S$ such that $u\nleq t$, and choose a maximal such $t$. Then $t\notin U$, hence $t\in\Grm(S)$. Let $\hat{t}$ be the cogerm of $t$ in~$S$. There is an integer $m\in \N$ and elements $t_i\in S$, for $0\leq i\leq m$,  such that
$$[\,t,\hat{t}\,]_S=\{\,t=t_0<t_1\ldots<t_m=\hat{t}\,\}\mpoint$$
If $t_i\notin U$ for some $i\in\{1,\ldots,n\}$, then $t_i\in \Grm(S)$, and $t_i>t$. By Lemma~\ref{lemme racines} $t_i>\hat{t}=t_n$, a contradiction. Thus $t_i\in U$ for any $i\in\{1,\ldots,n\}$. In particular $u\leq t_1$. If $u<t_1$, then $u\leq t$, since $t=t_0$ is the greatest element of $\,]\,.\,,t_1[_S$. This contradiction shows that $u=t_1$ is the smallest element of $]\hat{s},\,.\,[_U$.
Then $t$ is the greatest element of $\,]\,.,t_1[_S=\,]\,.,u[_S$. It follows that $t$ is unique, that is, there is a unique element $t$ of $]\hat{s},\,.\,[_S$ maximal subject to $u\nleq t$. In other words
$$]\hat{s},\,.\,[_S=]\hat{s},t]_S\msqcup[u,\,.\,[_S\mvirg$$
and $t<u=t_1$.\par
If $]\hat{s},t[_S\neq\emptyset$, let $t'$ be a maximal element of this poset. Then $t\notin U$, for otherwise $t'\geq u>t$, since $t'>\hat{s}$. So $t'\in\Grm{S}$, and $t'<t$. By Lemma~\ref{lemme racines}, it follows that $t'\leq \hat{t}'<t$, where $\hat{t}'$ is the cogerm of $t'$ in $S$. Hence $t'=\hat{t}'$ by maximality of $t'$. But then
$$]t',\,.\,[_S=\{t\}\msqcup[u,\,.\,[_S\mvirg$$
hence $\Inf_S]t',\,.\,[_S=t\neq t'$, contradicting $t'\in\Grm(S)$.\par
It follows that $]\hat{s},t[_S=\emptyset$, thus $]\hat{s},\,.\,[_S=\{t\}\msqcup[u,\,.\,[_S$. Then $\Inf_S]\hat{s},\,.\,[_S=t\neq s$, contradicting $s\in\Grm(S)$.\par
This shows finally that $u\leq t$ for any $t\in]\hat{s},\,.\,[_S$, thus $u\leq \hat{s}=\Inf_S]\hat{s},\,.\,[_S$. Hence
\begin{equation}\label{inf}
\hat{s}=\Inf_U]\hat{s},\,.\,[_U\mpoint
\end{equation}
Now $s$ is a germ of $U$, by (\ref{sup}),~(\ref{intervalle}), and ~(\ref{inf}). This completes the proof of Assertion~1.\mpn
For Assertion 2, let $r$ be a germ of $U$, and let $\hat{r}$ be its cogerm in $U$. Let $t\in S$ such that $t\leq x$ for any $x\in ]\hat{r},\,.\,[_S$. Then for any $u\in U$ with $u\leq t$, we have $u\leq v$ for any $v\in]\hat{r},\,.\,[_U$. Hence $u\leq\Inf_U]\hat{r},\,.\,[_U=\hat{r}$. Thus $U_{\leq t}\subseteq U_{\leq \hat{r}}$, thus $t\leq\hat{r}$ since $U$ detects $S$ by Lemma~\ref{detecte}. This shows that 
\begin{equation}\label{inf2}
\hat{r}=\Inf_S]\hat{r},\,.\,[_S\mpoint
\end{equation}
Now let $t\in S$ with $t\leq\hat{r}$. Then $U_{\leq t}\subseteq \,]\,.\,,\hat{r}]_U=\,]\,.\,,r[_U\msqcup [r,\hat{r}]_U$, and there are two cases: either $U_{\leq t}\subseteq \,]\,.\,,r[_U$, and then $U_{\leq t}\subset U_{\leq r}=\,]\,.\,,r]_U$, thus $t< r$ as $U$ detects $S$. Or $U_{\leq t}$ has a greatest element $u\in[r,\hat{r}]_U$. In this case $U_{\leq t}=\,]\,.\,,u]_U=U_{\leq u}$, thus $t=u$. In other words
\begin{equation}\label{sous hatr}
\,]\,.\,,\hat{r}]_S=\,]\,.\,,r[_S\msqcup[r,\hat{r}]_U\mvirg
\end{equation}
and in particular $[r,\hat{r}]_S=[r,\hat{r}]_U$.\par
Now let $t\in S$ with $t\geq r$. By Theorem~\ref{S dans G}, there are two cases: In the first case there exists $B\subseteq U$ such that $U_{\leq t}=U_{\leq B}$. Then in particular $r\in U_{\leq B}$, that is $B\subseteq [r,\,.\,[_U=[r,\hat{r}]_U\msqcup ]\hat{r},\,.\,[_U$. If $B\subseteq ]\hat{r},\,.\,[_U$, then $\hat{r}\in U_{\leq B}=U_{\leq t}$, thus $\hat{r}\leq t$. Otherwise the set $B$ has a smallest element $u\in [r,\hat{r}]_U$, and in this case $U_{\leq B}=U_{\leq u}=U_{\leq t}$, hence $t=u\in [r,\hat{r}]_U$.\par
The other case is when there exists a germ $r'$ of $U$ such that $U_{\leq t}=\,]\,.\,,r'[_U$. Then $r<r'$ since $r\in U_{\leq t}$, and $\hat{r}<r'$ by Lemma~\ref{lemme racines}. It follows that $U_{\leq\hat{r}}\subseteq U_{\leq t}$, so $\hat{r}\leq t$ in this case also. We get finally that
\begin{equation}\label{sur r}
[r,\,.\,[_S=[r,\hat{r}]_U\msqcup]\hat{r},\,.\,[_S\mpoint
\end{equation}
It follows from (\ref{inf2}), (\ref{sous hatr}) and (\ref{sur r}) that $r$ is a germ of $S$ if and only if $r=\Sup_S\,]\,.\,,r[_S$, and in this case $\hat{r}$ is the cogerm of $r$ in $S$. This is Case (a) of Assertion~2.\par
And if $r\neq\Sup_S\,]\,.\,,r[_S$, there exists $s\in S$ such that $s\geq x$ for any $x\in\,]\,.\,,r[_S$, but $s\ngeq r$. In particular $\,]\,.\,,r[_U\subseteq U_{\leq s}$.\par
If there exists a subset $B$ of $U$ such that $U_{\leq s}=U_{\leq B}$, then $\,]\,.\,,r[_U\subseteq U_{\leq b}$, for any $b\in B$, thus $b\geq\Sup_U\,]\,.\,,r[_U=r$. Hence $r\in U_{\leq B}$, that is $r\in U_{\leq s}$, contradicting $s\ngeq r$. \par
So there is a germ $r'$ of $U$ such that $U_{\leq s}=\,]\,.\,,r'[_U$. It follows that $\,]\,.\,,r[_U\subseteq \,]\,.\,,r'[_U\subseteq \,]\,.\,,r']_U$, so $r'\geq \Sup_U\,]\,.\,,r[_U=r$. If $r'>r$, then $r\in\,]\,.\,,r'[_U=U_{\leq s}$, so $r\leq s$, a contradiction. Hence $r'=r$, and $U_{\leq s}=\,]\,.\,,r[_U\subset U_{\leq r}$, so $s<r$. If $t\in\,]\,.\,,r[_S$, then $U_{\leq t}\subseteq \,]\,.\,,r[_U=U_{\leq s}$, so $t\leq s$. It follows that $s$ is the greatest element of $\,]\,.\,,r[_S$, so $\Sup_S\,]\,.\,,r[_S=s<r$, and $r$ is not a germ of~$S$. Moreover $s\notin U$ since $\Sup_U\,]\,.\,,r[_U=r$. Finally, the previous discussion shows that $s$ is the only element of $S$ such that $s\geq x$ for any $x\in\,]\,.\,,r[_S$, but $s\ngeq r$. In particular $r$ is the smallest element of $]s,\,.\,[_S$.  Now (\ref{inf2}), (\ref{sous hatr}) and (\ref{sur r}) show that $s$ is a germ of $S$, with cogerm $\hat{r}$ in $S$. This is case (b) of Assertion~2, and completes the proof of Theorem~\ref{extension racines}.  \endpf 
\begin{mth}{Corollary}\label{G donne U} \begin{enumerate} 
\item Let $T$ be a finite lattice. Set $U=T-\Grm(T)$, considered as a full subposet of $T$. Then $T\cong G(U)$.
\item Let $U$ be a finite poset, and let $\sou{U}$ denote its isomorphic image in its germ closure $T=G(U)$. Then $\sou{U}=T-\Grm(T)$.
\item Let $U$ and $V$ be finite posets. Let $\varphi:G(U)\to G(V)$ be an isomorphism of posets. Then $\varphi(\sou{U})=\sou{V}$, and in particular $U$ and $V$ are isomorphic.
\item In particular, the restriction $\varphi\mapsto \varphi_{\mid\sou{U}}$ induces a group isomorphism $\Aut\big(G(U)\big)\cong \Aut(U)$.
\end{enumerate}
\end{mth}
\pf For Assertion 1, the poset $T$ is clearly a germ extension of $U=T-\Grm(T)$. By Theorem~\ref{germ closure}, the map $j:t\in T\mapsto U_{\leq t}\in G(U)$ is a poset isomorphism from $T$ to a full subposet of $G(U)$ containing $\sou{U}$. All we have to show is that this map is surjective. Clearly for $B\subseteq U$, we have $j(\mathop{\meet}_{b\in B}\limits b)=U_{\leq B}$, so $\Lambda(U)$ is contained in the image of $j$. Now if $r\in \Grm(U)$, then $r\notin\Grm(T)$ (since $r\in U$), hence $\,]\,.\,,r[_T$ has a greatest element $s\in T-U$, by Theorem~\ref{extension racines}. Then $\,]\,.\,,r[_U=U_{\leq s}=j(s)$. Hence $\Blb(U)$ is contained in the image of $j$, completing the proof of Assertion~1.\mpn
For Assertion~2, by Theorem~\ref{germ closure}, the poset $T=G(U)$ is a germ extension of~$\sou{U}$. Hence $T=\sou{U}\mcup\Grm(T)$. So all we have to show is that $\sou{U}\mcap \Grm(T)=\emptyset$. Let $s\in\sou{U}\mcap \Grm(T)$. Then $s=\,]\,.\,,u]_U$, for some $u\in U$, and by Theorem~\ref{extension racines}, the element $u$ is a germ of $U$. This means that $\,]\,.\,,u[_U\in G(U)$. Since any lower-subset of $U$ properly contained in $\,]\,.\,,u]_U$ is contained in $\,]\,.\,,u[_U$, it follows that $\,]\,.\,,s[_T$ has a greatest element $t=\,]\,.\,,u[_U$. So $\Sup_T\,]\,.\,,s[_T=t<s$, so $s$ is not a germ of $T$, and Assertion 2 follows by contradiction. \mpn
Assertion 3 follows as well: if $\varphi:G(U)\to G(V)$ is an isomorphism of posets, then since $\sou{U}=G(U)-\Grm\big(G(U)\big)$ and $\sou{V}=G(V)-\Grm\big(G(V)\big)$, it follows that $\varphi(\sou{U})=\sou{V}$. Hence $U$ and $V$ are isomorphic.\par
In the case $U=V$, this shows that $\varphi(\sou{U})=\sou{U}$ for any automorphism $\varphi$ of the poset $G(U)$. By Theorem~\ref{germ closure}, the resulting group homomorphism $\Aut\big(G(U)\big)\to \Aut(\sou{U})\to \Aut(U)$ is injective. It is also surjective, since any automorphism $\alpha$ of $U$ extends to an automorphism of $G(U)$: indeed $\alpha(U_{\leq B})=U_{\leq\alpha(B)}$, for $B\subseteq U$, and moreover $\alpha\big(\,]\,.\,,r[_U\big)=\,]\,.\,,\alpha(r)[_U$ for any $r\in U$, and $\alpha\big(\Grm(U)\big)=\Grm(U)$. This completes the proof.\vspace{-2ex}
\findemo
\section{Germ extensible subsets of a lattice}
\begin{mth}{Definition and Notation}\label{def extensible} Let $T$ be a finite lattice. A full subposet~$U$ of $T$ is called {\em germ extensible in $T$} if the natural map
$$\nu :s\in G(U)\mapsto \join_{u\in s}u\in T\vspace{-2ex}$$
is injective. We denote by $\sur{G}(U)$ the image of $\nu$.
\end{mth}
\begin{mth}{Theorem}\label{germ extensible} Let $T$ be a finite lattice, and $U$ be a full subposet of $T$. Then $U$ is germ extensible in $T$ if and only if
\begin{equation}\label{extensible} \forall r\in\Grm(U),\;\;r>\joinb{u\in U}{u<r}u\;\;\hbox{in} \;\;T\mpoint
\end{equation}
In this case, the map $t\in\sur{G}(U)\mapsto \{u\in U\mid u\leq t\}$ is inverse to the bijection $\nu:G(U)\to\sur{G}(U)$.
\end{mth}
\pf If $U$ is germ extensible in $T$, let $r\in\Grm(U)$. Then $\,]\,.\,,r[_U$ and $\,]\,.\,,r]_U$ are distinct elements of $G(U)$. Thus $\nu\big(\,]\,.\,,r]_U\big)=r>\nu\big(\,]\,.\,,r[_U\big)=\joinb{u\in U}{u<r}u$. So Condition~\ref{extensible} is necessary.\par
Conversely, let $B\subseteq U$, and set $s=U_{\leq B}\in G(U)$ and $\meet B=\meet_{b\in B}b\in T$.  Then $\nu(s)=\join_{u\in U_{\leq B}}u\leq\meet B$, and for any $v\in U$
$$v\leq \meet B\implies v\in U_{\leq B}\implies v\leq \nu(U_{\leq B})=\nu(s)\implies v\leq \meet B\mpoint$$
Thus $s=\{v\in U\mid v\leq \nu(s)\}$.\par
Now let $r$ be a germ of $U$, set $s=\,]\,.\,,r[_U\in G(U)$. If $r>\joinb{u\in U}{u<r}u=\nu(s)$, then for any $v\in U$
$$v\in s \implies v<r\implies v\leq \joinb{u\in U}{u<r}u=\nu(s)\implies v<r\implies v\in s\mvirg$$
so $s=\{v\in U\mid v\leq\nu(s)\}$ in this case also.\findemo
\pagebreak[3]
\begin{mth}{Theorem} \label{partition}Let $T$ be a finite lattice, and let $S$ be a subset of $T$. Then there exists a unique germ extensible subset $U$ of $T$ such that $U\subseteq S\subseteq \sur{G}(U)$, namely
$$U=S-\{s\in\Grm(S)\mid s=\joinb{t\in S}{t<s}t\}\mpoint$$
In other words, the lattice of subsets of $T$ is the disjoint union of the intervals $[U,\sur{G}(U)]$, when $U$ runs through all germ extensible subsets of $T$.
\end{mth}
\pf Let $S\subseteq T$, and suppose that $U$ is a germ extensible subset of $T$ such that $U\subseteq S\subseteq \sur{G}(U)$. Since $S\subseteq \sur{G}(U)$, any element $t$ of $S$ is a join of elements of $U$, that is $t=\joinb{u\in U}{u\leq t}u$. It follows that for any $s\in S$
$$\joinb{t\in S}{t<s}t=\joinb{u\in U}{u<s}u\mpoint$$
Let $s\in \Grm(S)$.  If $s\in U$, then $s\in \Grm(U)$ by Theorem~\ref{extension racines}, hence $s>\joinb{u\in U}{u<s}u$ by Theorem~\ref{germ extensible}. In other words if $s=\joinb{t\in S}{t<s}t$, then $s\notin U$. Thus 
$$U\subseteq S-\{s\in\Grm(S)\mid s=\joinb{t\in S}{t<s}t\}\mpoint$$
Conversely, let $V=S-\{s\in\Grm(S)\mid s=\joinb{t\in S}{t<s}t\}$. Then $S$ is a germ extension of $V$, so by Theorem~\ref{germ closure}, there is an isomorphism of $S$ onto a full subposet $S'$ of $G(V)$ such that $\sou{V}\subseteq S'\subseteq G(V)$. \par
If $v\in \Grm(V)$, there are two cases. First if $v\in\Grm(S)$, then $v>\joinb{t\in S}{t<v}t$ by definition of $V$. Then {\em a fortiori} $v>\joinb{w\in V}{w<v}w$. Now if $v\notin \Grm(S)$, then $\sou{v}=\,]\,.\,,v]_V\notin\Grm(S')$, so $\,]\,.\,,\sou{v}[_{S'}$ has a greatest element $s'$ by Theorem~\ref{extension racines}. Thus $\,]\,.\,,v[_S$ has a greatest element $s$, and then $\joinb{w\in V}{w<v}w\leq s<v$ in this case also. \par
It follows from Theorem~\ref{germ extensible} that $V$ is a germ extensible subposet of $T$. The inclusions $\sou{V}\subseteq S'\subseteq G(V)$ now read $V\subseteq S\subseteq \sur{G}(V)$ in $T$, and we have $\rule{0ex}{2.5ex}U\subseteq V\subseteq S\subseteq \sur{G}(U)\mcap\sur{G}(V)\subseteq T$. \par
If $v\in V-U$, then $v\in\Grm(V)$, so $v>\joinb{w\in V}{w<v}w$. But on the other hand $v\in\sur{G}(U)$, so $v=\joinb{u\in U}{u\leq v}u=\joinb{u\in U}{u< v}u$. Hence $\joinb{u\in U}{u< v}u>\joinb{w\in V}{w<v}w$, a contradiction, since any element $w$ of $V$ is the join of the elements $u\leq w$ of $U$. It follows that $U=V$, as was to be shown.\endpf
\section{Germs and lattices}\label{germs and lattices}
\npar Let $T$ be a finite lattice, and let $E=\Irr(T)$ be the set of join-irreducible elements of $T$. The following constructions are introduced in~\cite{corfun-simplicity} (Notation 2.3 and 2.6). First, we denote by
$$\Lambda E=\{t\in T\mid t=\meet_{\substack{e\in E\\e\geq t}}e\}$$
the set of elements of $T$ which are equal to a meet of join-irreducible elements of $T$. Moreover, for in T, we set
$$r(t)=\join_{\substack{e\in E\\e<t}}e\;\;\hbox{and}\;\;\sigma(t)=\meet_{\substack{e\in E\\e>t}}e\mpoint$$
So $r(t)\leq t$, with equality if and only if $t\notin E$. And if $t\in E$, then $r(t)$ is the largest element of $]\,.\,,t[_T$. In particular, it follows from Theorem~\ref{germ extensible} that $E$ is a germ extensible subset of $T$. \par
On the other hand $\sigma(t)\geq t$ with equality and only if either $t\in \Lambda E-E$, or $t\in E$ and $t$ is equal to the meet of the elements of $E\,\cap\, ]t,\,.\,[_T$.\par
For any $t\in T$, we denote by $r^\infty(t)$ the limit of the decreasing sequence $t\geq r(t)\geq r^2(t)\geq \ldots$, and by $\sigma^\infty(t)$ the limit of the increasing sequence $t\leq\sigma(t)\leq\sigma^2(t)\leq\ldots$.\par
Finally (\cite{corfun-simplicity} Notation 2.10), we set
\begin{align}\label{def GT}
G^\sharp_T&=\{t\in T\mid t=r^\infty\sigma^\infty(t)\},\notag\\
\widehat{G}_T&=G^\sharp-\Lambda E,\\
G_T=E\sqcup G_T^\sharp&=E\sqcup (\Lambda E-E)\sqcup \widehat{G}_T=\Lambda E\sqcup \widehat{G}_T\mpoint\notag
\end{align}
We just saw that $E$ is a germ extensible subset of $T$. The following theorem shows that the corresponding set $\sur{G}(E)$ is equal to the set $G_T$ introduced above. In particular, it only depends on the poset $E$.
\begin{mth}{Theorem} \label{claim}Let $T$ be a finite lattice, and $E$ the full subposet of join-irreducible elements of $T$. Then $E$ is germ extensible in $T$, and moreover $\sur{G}(E)=G_T$.
\end{mth} 
\pf For $t\in T$, we set $\alpha(t)=\{e\in E\mid e\leq t\}$. \mp
We have $G_T=\Lambda E\sqcup \widehat{G}_T$. So if $t\in G_T$, then either $t\in\Lambda E$, that is $t=\meetb{e\in E}{e\geq t}e$, and $\alpha(t)=E_{\leq B}$, where $B=\{e\in E\mid e\geq t\}$. 
Hence $\alpha(t)\in G(E)$ in this case. \par
Otherwise $t\in\widehat{G}_T$. Then $t\notin \Lambda E$, so we have a sequence
$$t<\sigma(t)<\sigma^2(t)<\ldots<\sigma^n(t)=\sigma^\infty(t)\mpoint$$
All the terms different from $t$ of this sequence are in $E$: indeed, if $\sigma^i(t)\notin E$ for $i\geq 1$, then $r\sigma^i(t)=\sigma^i(t)$, so 
$$t=r^\infty\sigma^\infty(t)\geq r^\infty\sigma^i(t)=\sigma^i(t)\mvirg$$
contradicting $t<\sigma^i(t)$. So $\sigma^{i-1}(t)\leq r\sigma^i(t)<\sigma^i(t)$ for $i\geq 1$. Moreover if $\sigma^{i-1}(t)< r\sigma^i(t)$, then there are two cases: either $r\sigma^i(t)\notin E$, and then
$$t\leq \sigma^{i-1}(t)<r\sigma^i(t)=r^\infty\sigma^i(t)\leq r^\infty\sigma^\infty(t)=t\mvirg$$
a contradiction. Or $r\sigma^i(t)\in E$, and then $r\sigma^i(t)\geq\sigma\big(\sigma^{i-1}(t)\big)=\sigma^i(t)$, contradicting $r\sigma^i(t)<\sigma^i(t)$. Hence $r\sigma^i(t)=\sigma^{i-1}(t)$ for any $i\in\{1,\ldots,n\}$.\mp

We set $\gamma=\sigma(t)>t$. Then $\gamma\in E$, and $\gamma$ is a germ of $E$. Indeed:
\begin{itemize}
\item The element $t$ is the greatest element of $\,]\,.\,,\gamma[_T$, so 
$$\,]\,.\,,\gamma[_E=\{f\in E\mid f<\gamma\}=\{f\in E\mid f\leq t\}=\alpha(t)\mpoint$$
Hence if $f\in E$ and $f\geq g$ for any $g\in\,]\,.\,,\gamma[_E$, then $f\geq \joinb{e\in E}{e\leq t}e=t$, hence $f>t$ and $f\geq \gamma$. Thus $\gamma=\Sup_E\,]\,.\,,\gamma[_E$.
\item Let $e_i=\sigma^i(t)$, for $i\geq 1$. If $e\in E$ and $e\geq \gamma=e_1$,  and either $e>e_n$, or there exists a largest integer $i\in\{1,\ldots,n\}$ such that $e_i\leq e$. If $e_i<e$, then $\sigma(e_i)=e_{i+1}\leq e$, contradicting the definition of $i$. Hence $e_i=e$ for some $i\in\{1,\ldots,n\}$, or $e>e_n$. Similarly, since $e_i=r(e_{i+1})$ for $i\in \{0,\ldots,n-1\}$, if $e\in E$ and $e\leq e_n$, then $e$ is equal to $e_i$ for some $i\in\{1,\ldots,n\}$, or $e<\gamma$.
\item Finally $e_n=\sigma(e_n)=\meetb{e\in E}{e>e_n}e$. Thus if $f\in E$ and $f\leq e$ for all $e\in]e_n,\,.\,[_E$, then $f\leq e_n$. In other words $e_n=\Inf_E]e_n,\,.\,[_E$.
\end{itemize}
This shows that $\gamma$ is a germ of $E$, with cogerm $e_n$. Now $\alpha(t)=\,]\,.\,,\gamma[_E$, so $\alpha(t)\in G(E)$ in this case also. It follows that $\alpha(t)\in G(E)$ for any $t\in G_T$. \mp
Conversely, for $s\in G(E)$, set $\nu(s)=\join_{e\in s}e$. Then $\nu(s)\in G_T$: indeed, either $s=E_{\leq B}$ for some $B\subseteq E$, and then $\nu(s)=\meet_{b\in B}b\in\Lambda E$. Or $s=\,]\,.\,,g[_E$, for $g\in \Grm(E)$. In this case $t=\nu(s)=r(g)\in\widehat{G}_T$: indeed, if $\hat{g}$ is the cogerm of $g$ in $E$, and if $[g,\hat{g}]_E=\{g=e_0<e_1<\ldots<e_n=\hat{g}\}$, then clearly $e_i=\sigma^{i+1}(t)$ for $0\leq i\leq n$. Moreover $\sigma(\hat{g})=\hat{g}$, since $\sigma(\hat{g})=\meetb{f\in E}{f>\hat{g}}f$ is equal to the join of all elements $e\in E$ such that $e\leq f$ for all $f\in]\hat{g},\,.\,[_E$, and $\Inf_E]\hat{g},\,.\,[_E=\hat{g}$. Then clearly again $r(e_i)=e_{i-1}$ for $1\leq i\leq n$, and $t=r(g)$.\par
So we have the maps $\alpha$ and $\nu$
$$\xymatrix{
G_T\ar[r]<.4ex>^-\alpha&G(E)\ar[l]<.4ex>^-\nu\mvirg
}
$$
which are obviously maps of posets. Moreover $\nu\circ\alpha(t)=t$ for any $t\in G_T$, since $t=\joinb{e\in E}{e\leq t}e$ for any $t\in T$. Similarly, if $s\in G(E)$, then $\alpha\circ \nu(s)=s$: indeed either $s=E_{\leq B}$ for some $B\leq U$, and then $\nu(s)=\meet_{b\in B}b$. In this case $e\leq \nu(s)$ for $e\in E$ if and only if $e\in E_{\leq B}=s$, so $\alpha\circ \nu(s)=s$. Or $s=\,]\,.\,,g[_E$ for some $g\in \Grm(E)$, and then $\nu(s)=r(g)$, so that $e\leq \nu(s)$ for $e\in E$ if and only if $e<g$, that is $\alpha\circ \nu(s)=s$ in this case also. This completes the proof of Theorem~\ref{claim}, since by Definition~\ref{def extensible}, the image of $\nu$ is equal to $\sur{G}(E)$.\endpf
\pagebreak[3]
\npar Let $U$ be a finite poset, and $\Idown(U)$ be the lattice of lower-subsets of~$U$. The irreducible elements of $\Idown(U)$ are the subsets of the form $]\,.\,,u]_U$, for $u\in U$, so we can identify $\Irr\big(\Idown(U)\big)$ with $U$. We will now show that $G(U)$ is equal to the subset $G_{\Idown(U)}$ of the lattice $\Idown(U)$ introduced in (\ref{def GT}).\par
First let $B$ be a subset of $U$. Then $\mathop{\bigcap}_{b\in B}\limits]\,.\,,b]_U=U_{\leq B}$. In other words the set $\Lambda U$ of intersections of irreducible elements of $\Idown(U)$ is equal to the set $\Lambda(U)$ of Definition~\ref{def G}.\par
Now let $S\in\widehat{G}_{\Idown(U)}$. Then $S$ is not a meet of join-irreducible elements of $\Idown(U)$, so $S\neq U_{\leq B}$, for any $B\subseteq U$. Now in the proof of Theorem~\ref{claim}, we saw that in the sequence
$$S<\sigma(S)<\sigma^2(S)<\ldots<\sigma^n(S)=\sigma^\infty(S)\mvirg$$
all the terms different from $S$ are join-irreducible in $\Idown(U)$. Moreover $r\sigma^i(S)=\sigma^{i-1}(S)$ for $i\geq 1$.\par
It follows that there is a sequence $u_1,\ldots,u_n$ of elements of $U$ such that $\sigma^i(U)=]\,.\,,u_i]_U$, for $i\geq 1$. In particular $u_1<u_2<\ldots<u_n$. Moreover since $r\big(]\,.\,,u]_U\big)=]\,.\,,u[_U$ for any $u\in U$, we have $S=]\,.\,,u_1[_U$, and $]\,.\,,u_i[_U=]\,.\,,u_{i-1}]_U$ for $i\geq 2$. In other words $S=]\,.\,,u_1[_U$, and $u_{i-1}$ is the largest element of $]\,.\,,u_i[_U$, for $i\geq 2$.\par
Now $\sigma(S)=]\,.\,,u_1]_U=\mathop{\bigcap}_{\substack{u\in U\\S\subset]\,.\,,u]_U}}\limits]\,.\,,u]_U$. Since $S$ is not irreducible in $\Idown(U)$, saying that $S\subset\,]\,.\,,u]_U$ is equivalent to saying that $S\subseteq\,]\,.\,,u]_U$, that is, $u$ is an upper bound of $S$. Thus $u_1$ is the smallest upper bound of $S$, i.e. $u_1=\Sup_U]\,.\,,u_1[_U$. Similarly for $1\leq i<n$, saying that $\sigma\big(]\,.\,,u_i]_U\big)=]\,.\,,u_{i+1}]_U$ amounts to saying that $u_{i+1}$ is the smallest element of $]u_i,\,.\,[_U$. Since moreover $r\sigma^i(S)=\sigma^{i-1}(S)$ for $i\geq 1$, it follows that $]\,.\,,u_i[_U=]\,.\,,u_{i-1}]_U$ for $i\geq 2$, i.e. $u_{i-1}$ is the largest element of $]\,.\,,u_i]_U$, for $i\geq 2$. Finally saying that $]\,.\,,u_n]_U=\sigma\big(]\,.\,,u_n]_U\big)=\mathop{\bigcap}_{\substack{u\in U\\u>u_n}}\limits ]\,.\,,u]_U$ amounts to saying that $u_n$ is the greatest lower bound of $]u_n,\,.\,[_U$, i.e. $u_n=\Inf_U]u_n,\,.\,[_U$.\par
This discussion shows that $u_1$ is a germ of $U$, with cogerm $u_n$. Hence $S=]\,.\,,u_1[_U$ belongs to the set $\widehat{G}(U)$ of Definition~\ref{def G}. Conversely, if $u_1$ is a germ of $U$, with cogerm $u_n$, it is straightforward to reverse the above arguments and check that $S=]\,.\,,u_1[_U$ belongs to $\widehat{G}_{\Idown(U)}$, as defined in~\ref{def GT}, with $\sigma^\infty(S)=]\,.\,,u_n]_U$. This shows that Notation~\ref{def G} is consistent with Notation~\ref{def GT}.
\section{Examples}\label{examples}
\npar The empty subset of a finite (non empty) lattice $T$ is germ extensible (by Theorem~\ref{germ extensible}, since the emptyset has no germs at all). Moreover $\sur{G}(\emptyset)=\{0\}$, where $0$ is the smallest element of $T$.
\npar Let $U=\{u_1<u_2<\ldots<u_n\}$ be a totally ordered poset of cardinality $n>0$. Then the only germ of $U$ is $u_1$, with cogerm $u_n$. The poset $G(U)$ is equal to $\{V_0\subset V_1\subset \ldots\subset V_n\}$, where $V_i=\{u_1,\ldots,u_i\}$ for $0\leq i\leq n$ (so $V_0=\emptyset$). If $U$ is a subposet of a finite lattice $T$, then $U$ is germ extensible in $T$ if and only if $u_1$ is not equal to the smallest element $0$ of $T$. In this case $\sur{G}(U)=\{0\}\sqcup U$.
\npar Let $U=\{u_1,\ldots,u_n\}$ be a discrete poset of cardinality $n\geq 2$. Then $U$ has no germs, and $G(U)=\{\emptyset\}\sqcup \big\{\{u\}\mid u\in U\big\}\sqcup\{U\}$. If $U$ is a full subposet of a finite lattice $T$, then $U$ is germ extensible in $T$, and $\sur{G}(U)=\{0\}\sqcup U\sqcup\{v\}$, where where $0$ is the smallest element of $T$ and $v=u_1\join u_2\join\ldots\join u_n$ in $T$.
\npar Let $U=\{a<c>b\}$ be a connected poset of cardinality $3$ with two minimal elements. Then $\Grm(U)\!=\!\{c\}$, and $G(U)\!=\!\big\{\emptyset,\{a\},\{b\},\{a,b\},\{a,b,c\}\big\}$. If $U$ is a full subposet of a finite lattice $T$, then $U$ is germ extensible in $T$ if and only if $c>a\join b$ in $T$. In this case $\sur{G}(U)=\{0,a,b,a\join b,c\}$.
\npar Let $V=U\op=\{a>c<b\}$ be the opposite poset of the previous example. Then $\Grm(V)\!=\!\{c\}$, and $G(V)\!=\!\big\{\emptyset,\{c\},\{b,c\},\{a,c\},\{a,b,c\}\big\}$. If $V$ is a full subposet of a finite lattice $T$, then $V$ is germ extensible in $T$ if and only if $c$ is not equal to the smallest element $0$ of $T$. In this case $\sur{G}(V)=\{0,c,a,b,a\join b\}$.
\npar In the previous two examples, there is an isomorphism of posets $G(U\op)\cong G(U)\op$. This is a general phenomenon, that will be established in a joint forthcoming paper with Jacques Thévenaz (\cite{corfun-duality}).
\npar In all the previous examples, the poset of join-irreducible elements of the lattice $G(U)$ is isomorphic to $U$, but this need not be true in general. For example
$$\hbox{if}\;U=\vcenter
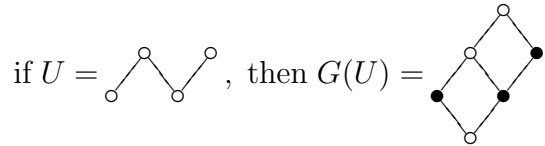
{\xymatrix@C=1.8ex@R=1.5ex@M=0.3ex{&\pointcreux\ar@{-}[ld]\ar@{-}[rd]&&\pointcreux\ar@{-}[ld]\\
\pointcreux&&\pointcreux&
}}\;,\;\hbox{then}\; 
G(U)=\vcenter{\xymatrix@C=1.8ex@R=1.5ex@M=0.3ex{&&\pointcreux\ar@{-}[ld]\ar@{-}[rd]&\\
&\pointcreux\ar@{-}[ld]\ar@{-}[rd]&&\pointplein\ar@{-}[ld]\\
\pointplein\ar@{-}[rd]&&\pointplein\ar@{-}[ld]\\
&\pointcreux
}}
$$
where the irreducible elements of $G(U)$ are the black ones.
\npar An example of a lattice $T$, its germs ($\,\pointcarre\,$), the poset $U$ of its irreducible elements ($\,\pointplein\,$), and the germ closure $G(U)$ with its germs ($\,\pointcarre\,$):
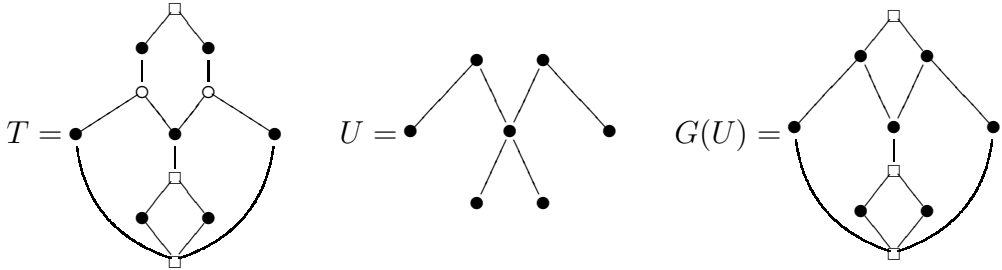
$$
T=\vcenter{\xymatrix@C=1.8ex@R=1.5ex@M=0.3ex{
&&&\pointcarre\ar@{-}[rd]\ar@{-}[ld]&&&\\
&&\pointplein\ar@{-}[d]&&\pointplein\ar@{-}[d]&&\\
&&\pointcreux\ar@{-}[rd]\ar@{-}[lld]&&\pointcreux\ar@{-}[rrd]\ar@{-}[ld]&&\\
\pointplein\ar@{-}@/_2ex/[rrrddd]&&&\pointplein\ar@{-}[d]&&&\pointplein\ar@{-}@/^2ex/[lllddd]\\
&&&\pointcarre\ar@{-}[rd]\ar@{-}[ld]&&&\\
&&\pointplein\ar@{-}[rd]&&\pointplein\ar@{-}[ld]&&\\
&&&\pointcarre&&&\\
}}
\hspace{.8cm}U=\vcenter{\xymatrix@C=1.8ex@R=1.5ex@M=0.3ex{
&&&&&&\\
&&\pointplein\ar@{-}[rdd]\ar@{-}[lldd]&&\pointplein\ar@{-}[rrdd]\ar@{-}[ldd]\\
&&&&&&\\
\pointplein&&&\pointplein\ar@{-}[ldd]\ar@{-}[rdd]&&&\pointplein\\
&&&&&&\\
&&\pointplein&&\pointplein\\
&&&&&&
}}
\hspace{.8cm}G(U)=\vcenter{\xymatrix@C=1.8ex@R=1.5ex@M=0.3ex{
&&&\pointcarre\ar@{-}[rd]\ar@{-}[ld]&&&\\
&&\pointplein\ar@{-}[rdd]\ar@{-}[lldd]&&\pointplein\ar@{-}[rrdd]\ar@{-}[ldd]\\
&&&&&&\\
\pointplein\ar@{-}@/_2ex/[rrrddd]&&&\pointplein\ar@{-}[d]&&&\pointplein\ar@{-}@/^2ex/[lllddd]\\
&&&\pointcarre\ar@{-}[rd]\ar@{-}[ld]&&&\\
&&\pointplein\ar@{-}[rd]&&\pointplein\ar@{-}[ld]\\
&&&\pointcarre&&&
}}
$$
\npar Let $T$ be a finite lattice, and let $s\in\Grm(T)$. Then $s=\Sup_T\,]\,.\,,s[_T$, hence $s=\joinb{t\in T}{t<s}t$. By Theorem~\ref{partition}, it follows that the set $U=T-\Grm(T)$ is germ extensible in~$T$, and such that $\sur{G}(U)=T$. Moreover this set $U$ is the only germ extensible subset of $T$ with this property. This yields in particular another proof of Assertion~1 of Corollary~\ref{G donne U}.\mpn
{\bf Acknowledgements:} This paper would not exist without my long term collaboration with Jacques Th\'evenaz on correspondence functors, nor his numerous comments and pertinent suggestions. ``{\em Merci Jacques !}''.


\vspace{3ex}
\begin{flushleft}
Serge Bouc\\
LAMFA-CNRS UMR 7352\\
Universit\'e de Picardie-Jules Verne\\
33, rue St Leu, 80039 Amiens Cedex 01\\
France\\
{\small\tt serge.bouc@u-picardie.fr}\\
{\small\tt http://www.lamfa.u-picardie.fr/bouc/}
\end{flushleft}

\end{document}